# How Archimedes showed that π is approximately equal to 22/7

Damini D.B.[1] and Abhishek Dhar[2*]

[1]*National Academy for Learning, Basaveshwaranagar, Bengaluru 560079*

[2]*International Centre for Theoretical Sciences, TIFR, Bengaluru 560089*

**Abstract**

The ratio of the circumference (*C*) of a circle to its diameter (*D*) is a constant number denoted by $\pi$ and is independent of the size of the circle. It is known that $\pi$ is an irrational number and therefore cannot be expressed as a common fraction. Its value is approximately equal to 3.141592. Since Archimedes was one of the first persons to suggest a rational approximation of 22/7 for $\pi$, it is sometimes referred to as Archimedes' constant. In this article, we discuss how Archimedes came up with his formula. Archimedes in fact proved that $223/71 < \pi < 22/7$. Here we provide an improved lower bound.

______________________________________

## 1. Background

The mathematical constant $\pi$ is defined as

$$\pi = \frac{C}{D}, \qquad (1)$$

where $C$ is the circumference and $D$ is the diameter of the circle. The value of *C/D* is constant, regardless of the size of the circle. It is known that $\pi$ is an irrational number and therefore cannot be expressed as a ratio of two integers. Its value correct to the sixth decimal place is 3.141592 and fractions such as 22/7 and 355/113 are often used to approximate it. In 2019, Emma Haruka Iwao from japan numerically computed about 31 trillion digits of π, and the record for memorizing the maximum number of digits (71,000) is held by Rajveer Meena from India. Around 200 BCE Archimedes found a simple method for estimating the value of π. He noted that a regular polygon circumscribed around a circle would have a perimeter larger than the circumference of the circle, while the one inside the circle would have a smaller perimeter. He then observed that as one increased the number of sides of the polygon, it would get closer to the circle. Finally, he used Pythagoras' theorem to find the perimeter of the polygons and got upper and lower bounds for the value of π. In this article, we discuss the basic ideas behind Archimedes' proof. The details of his original method can be found in [1].

In Sections (2, 3), we first briefly discuss the history of the computation of $\pi$ and some formulas that mathematicians have discovered to represent this number. A detailed history of $\pi$ can be found in [1,2]. Then in Section (4), we describe the steps of Archimedes' derivation and make some concluding remarks in Section (5).

______________________________________

[2]Corresponding author;  *Email: abhishek.dhar@icts.res.in

## 2. A brief history of the number π

The value of π was first calculated 4000 years ago. The ancient Babylonians and the Egyptians calculated an approximate value of π by actual physical measurements of the circumference or the area of a circle and they estimated that π had a value close to 3.

About 1500 years later, the Greek mathematician Archimedes first used mathematics to compute π and showed that it had a value between 22/7 and 223/71. He did this by determining the perimeter of a polygon inscribed within a circle (which is less than the circle's circumference), and that of another polygon outside the circle (which is greater than the circle's circumference). He then found the average perimeter of the two. Using a hexagon, a 12-sided polygon, a 24-sided polygon, a 48-sided polygon and a 96-sided polygon, his final estimate gave a range between 3.1408 and 3.1428.

During the 5$^{th}$ century CE, Indian mathematicians calculated a value of π that was accurate up to 5 decimal digits. Zu Chongzhi, a Chinese mathematician and astronomer, calculated an approximate value of π (355/113) using a 24576-gon about 700 years after Archimedes. His estimated value, 355/113, is approximately equal to 3.14159292. The Greek letter π was first used by mathematicians in the 1700s and was introduced by William Jones in 1706. It was derived from the first letter of the Greek word 'perimetros,' meaning circumference.

## 3. Formulas for $\pi$

The basic formula for $\pi = C/D$, which is obtained by dividing the circumference of the circle *C* by its diameter *D*. However, the actual physical measurement of the circumference and the diameter of a circle would involve a lot of errors, hence the determination of π using this formula is not the most reliable method. There are many different formulas that mathematicians have discovered for representing the irrational number π. Here we give some examples [1-4].

A. Infinite series representations -

(a) *The Madhava-Gregory-Leibniz series (~1300-1700 CE):*

$$\frac{\pi}{4} = 1 - \frac{1}{3} + \frac{1}{5} - \frac{1}{7} + \frac{1}{9} - \frac{1}{11} \cdots$$

(b) *Nilakantha series (~1500 CE):*

$$\frac{\pi}{4} = \frac{3}{4} + \frac{1}{2\times3\times4} - \frac{1}{4\times5\times6} + \frac{1}{6\times7\times8} - \frac{1}{8\times9\times10} \cdots$$

The estimate for π gets more accurate as we include more terms in the above series.

B. Continued fractions representation – The following form was obtained by Brouncker (~1660 CE):

$$\frac{4}{\pi} = 1 + \cfrac{1^2}{2+\cfrac{3^2}{2+\cfrac{5^2}{2+\cfrac{7^2}{2+\cdots}}}}$$

C. Infinite product representation - The following product form was proposed by Wallis (1656 CE):

$$\frac{\pi}{2} = \left(\frac{2}{1} \cdot \frac{2}{3}\right) \cdot \left(\frac{4}{3} \cdot \frac{4}{5}\right) \cdot \left(\frac{6}{5} \cdot \frac{6}{7}\right) \cdot \left(\frac{8}{7} \cdot \frac{8}{9}\right) \cdots$$

D. Infinite nested radical representation - This involves taking an infinite number of square roots and was first proposed by Viete (1593 CE):

$$\frac{2}{\pi} = \frac{\sqrt{2}}{2} \times \frac{\sqrt{2+\sqrt{2}}}{2} \times \frac{\sqrt{2+\sqrt{2+\sqrt{2}}}}{2} \cdots$$

As we will see, Archimedes' method leads to a similar representation.

## 4. Archimedes' method for finding an approximate value of $\pi$

As discussed above, the basic idea of Archimedes is quite simple. We can draw a polygon with $n$ sides that circumscribes the circle. Let the perimeter of such a polygon be denoted by $C_n$. We also inscribe a polygon inside the circle. Let the perimeter of such a polygon be denoted by $c_n$. Then, if $C$ denotes the circumference of the circle, it is clear that we would have

$$c_n < C < C_n \tag{2}$$

As we increase $n$, both $c_n$ and $C_n$ should get closer to $C$. We illustrate this in Fig. (1) for the cases with $n = 3, 6, 12$. It is clear that the polygon with 12 sides (dodecagon) has a circumference that is the closest to the circle.

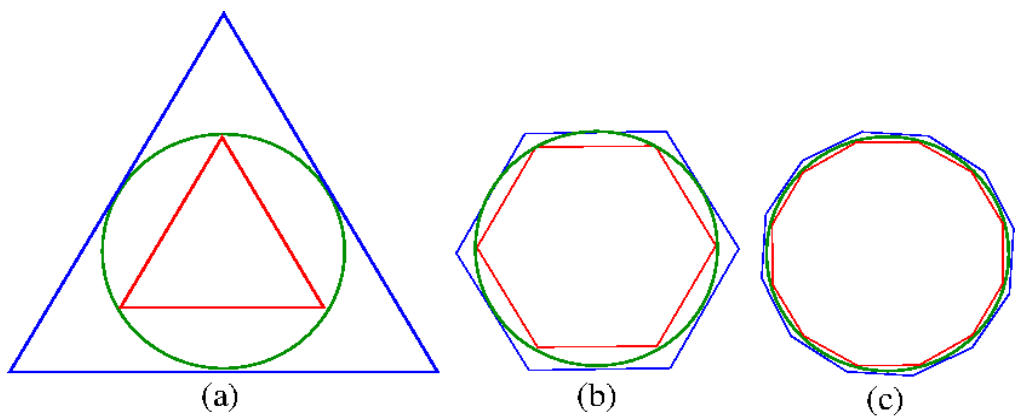

**Figure 1**: Here we show polygons of $n$ sides inscribed (red) and circumscribed (blue) around a circle. (a) Triangle ($n = 3$), (b) hexagon ($n = 6$), (c) dodecagon ($n = 12$).

Therefore, if we can find a way to compute the circumference of an $n$-sided polygon, we can arrive at an accurate value of $\pi$ simply by taking a large value for $n$. We now need a method to find the circumference of an *n*-sided polygon. We can do this using some trigonometry and the Pythagoras theorem though we must be careful not to use any result that already uses the value of $\pi$. From Fig. (2a,2b) it is clear that both the inscribed and circumscribed polygons for a circle of diameter 1 are respectively given by

$$c_n = n \sin\left(\frac{180}{n}\right), \quad C_n = n \tan\left(\frac{180}{n}\right) \qquad (3)$$

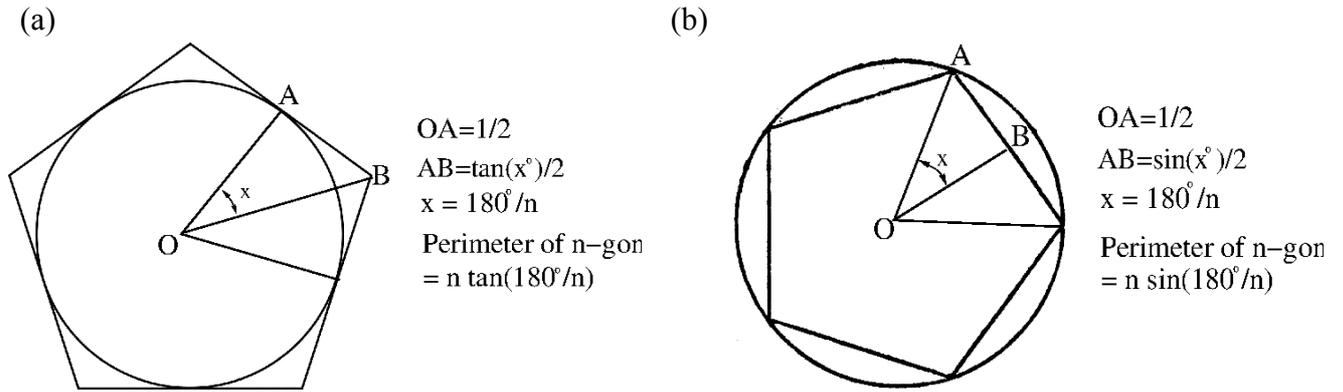

**Figure 2**: Here we show the calculation of (a) $C_n$ and (b) $c_n$, for *n*=5.

Now it turns out that for some particular choices of n, it is easy to evaluate these perimeters. For example, let us take *n=3*. Using Pythagoras theorem one can easily prove that sin(60°)=√3/2 so we get $c_3 = 3\sqrt{3}/2$. Next, using the formulas

$$\cos(x) = \sqrt{1 - sin^2 x}, \quad \tan(x) = \sin(x)/\cos(x) \qquad (4)$$

we get cos(60°) = 1/2 and tan(60°)= √3. Hence, we get $C_3 = 3\sqrt{3}$. If we now increase the number of sides from 3 by a factor of 2 to get a polygon with *n*=6, we now need the value of sin(30°) and tan(30°) to find the perimeters of the inscribed and circumscribed hexagons. Now cos(30°) is easy to find from cos(60°) using the trigonometric identity

$$\cos x = \sqrt{(1 + \cos 2x)/2}. \qquad (5)$$

This then gives us cos(30°)=√3/2 and sin(30°)=1/2. And then we go to *n*=12 for which we use the same rule to find cos(15°)=√(2 + √3)/2 and sin(15°)= √(2 − √3)/2. So the idea is to keep increasing *n* by powers of two and each time we can find the cosine, sine and tangent of the required angle by using the formulas in Eqs.(4,5). The results for *n=3,6,12,24,48,96* are given in Table-I. We see a pattern in the expressions and after sometime we can guess the next entry in the table. We also plot the numerical values (up to the 8[th] decimal place) and see that they become closer to $\pi$.

| n | $c_n = n \sin\left[\dfrac{180}{n}\right]$ | Numerical value (Approx.) | $C_n = n \tan\left[\dfrac{180}{n}\right]$ | Numerical value (Approx.) |
|---|---|---|---|---|
| 3 | $\dfrac{3\sqrt{3}}{2}$ | 2.59807621 | $3\sqrt{3}$ | 5.19615242 |
| 6 | $3$ | 3 | $2\sqrt{3}$ | 3.46410161 |
| 12 | $12\left[\dfrac{\sqrt{2-\sqrt{3}}}{2}\right]$ | 3.10582854 | $12\left[\dfrac{\sqrt{2-\sqrt{3}}}{\sqrt{2+\sqrt{3}}}\right]$ | 3.21539031 |
| 24 | $24\left[\dfrac{\sqrt{2-\sqrt{2+\sqrt{3}}}}{2}\right]$ | 3.13262861 | $24\left[\dfrac{\sqrt{2-\sqrt{2+\sqrt{3}}}}{\sqrt{2+\sqrt{2+\sqrt{3}}}}\right]$ | 3.15965994 |
| 48 | $48\left[\dfrac{\sqrt{2-\sqrt{2+\sqrt{2+\sqrt{3}}}}}{2}\right]$ | 3.13935020 | $48\left[\dfrac{\sqrt{2-\sqrt{2+\sqrt{2+\sqrt{3}}}}}{\sqrt{2+\sqrt{2+\sqrt{2+\sqrt{3}}}}}\right]$ | 3.14608622 |
| 96 | $96\left[\dfrac{\sqrt{2-\sqrt{2+\sqrt{2+\sqrt{2+\sqrt{3}}}}}}{2}\right]$ | 3.14103195 | $96\left[\dfrac{\sqrt{2-\sqrt{2+\sqrt{2+\sqrt{2+\sqrt{3}}}}}}{\sqrt{2+\sqrt{2+\sqrt{2+\sqrt{2+\sqrt{3}}}}}}\right]$ | 3.14271460 |

**Table 1**: This table shows the values of $c_n$ and $C_n$ for $n = 3, 6, 12, 24, 48, 96$

In Fig. (3a) we show a plot, which describes how $c_n$ and $C_n$ for a circle with unit diameter keep changing as we make *n* larger. We plot in Fig. (3b) a zoomed in version of Fig. (3a), the values 22/7=220/70 and 223/71 that were given by Archimedes and also the true value of $\pi$. We see that 22/7 is just a bit larger than $C_{96}$ while 223/71 is just a bit smaller than $c_{96}$.

It is clear that to get an accurate value of $\pi$, we have to take a very large value of *n*. In Fig. (3) we see that even with a 96-gon the difference is in the third decimal place. Up to 12 decimal places it is known that $\pi$ = 3.141592653589, while from our formula we find that $C_{24576}$ = 3.141592670702 and $c_{24576}$ = 3.141592645034. The Chinese mathematician Zu Chongzhi in fact used these results for the 24576-gon to guess the approximate rational value 355/113 which is about 3.141592920354.

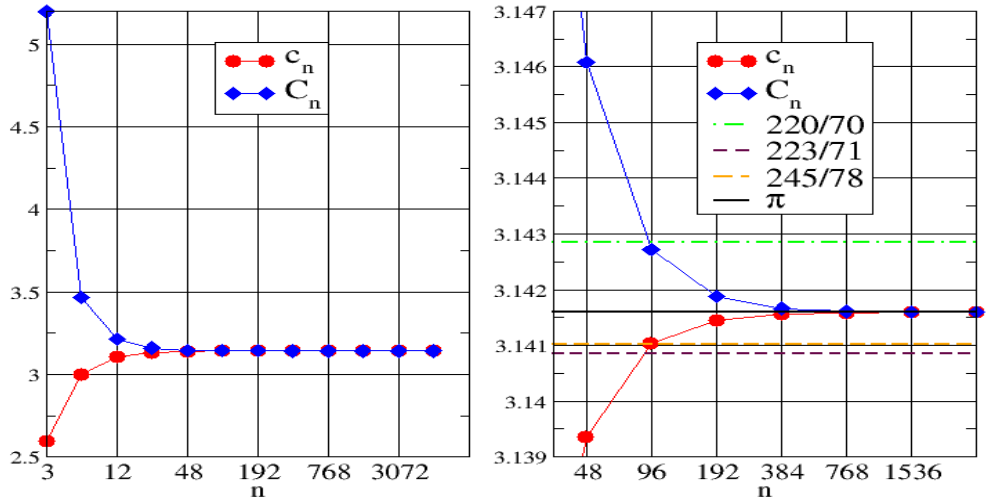

**Figure 3**: The numerical values obtained for $c_n$ and $C_n$ have been plotted as a function of $n$. The right panel shows a close-up so that we can see the values 22/7 and 223/71 obtained by Archimedes and also the true value of $\pi$. The improved lower bound 245/78 is also plotted.

## 5. Rational approximations

We see from Table 1 that that the perimeter of the polygons typically involves taking multiple square roots and is an irrational number. An interesting question is to find a rational number that is very close to this irrational number. Consider the value $c_{96} \approx 3.14103195$. A rational number approximating this (called a rational approximant) is 314103/100000. But can we get a rational number that has a smaller denominator? Archimedes obtained the value $223/71 \approx 3.140845$ which involves a smaller denominator and is still a good approximation. We are not sure how Archimedes made this guess.

However, one approach to find a good rational approximation for any number is to write it in the form of a continued fraction expansion, that we truncate after some term. We discuss this idea briefly and apply it to our problem.

A continued fraction expansion is given by the following form

$$a_0 + \cfrac{1}{a_1 + \cfrac{1}{a_2 + \cfrac{1}{a_3 + \cfrac{1}{a_4 + \cdots}}}}$$

where $a_0, a_1, \cdots$ are positive integers. For an irrational number, the series goes on for ever while for a rational number it would end after a finite number of terms. If we specify a number by its numerical value (up to some decimal places) there is a simple procedure to write down its continued fraction expansion. For illustration let us consider the number 3.14=157/50. We see that $a_0 = 3$. The remaining part of the expansion has the value 7/50. To find $a_1$, we have to take the integer part of 50/7. This gives us $a_1 = 7$. Repeating the procedure with the remainder, that is 1/7, we find $a_2 = 7$. Hence we can write

$$3.14 = 3 + \cfrac{1}{7 + \cfrac{1}{7}}.$$

So we get the series of rational approximants for 3.14 as: 3/1, 22/7, 157/50. If we perform this procedure for $c_{96} \approx 3.14103195$, we get

$$\frac{314103195}{100000000} = 3 + \cfrac{1}{7 + \cfrac{1}{11 + \cfrac{1}{25 + \cfrac{1}{1 + \cfrac{1}{25 + \cfrac{1}{1 + \cfrac{1}{27 + \cfrac{1}{13}}}}}}}}$$

Truncating the above series at various orders, we get the rational approximants 3/1, 22/7, 245/78, 6147/1957, and so on. This does not include Archimedes' value 223/71. However, notice that 245/78 ≈3.141025641 is in fact a better approximation to $c_{96}$ than 223/71. Moreover we see that 245/78 < $c_{96}$ and so we can use it as a lower bound. Applying the same procedure to $C_{96}$ gives us the approximants 3/1, 22/7, 3149/1002 and so on. In this case we see that 22/7 > $C_{96}$ and so we can use this value as an upper bound. Hence we finally get the somewhat improved bound

$$245/78 < \pi < 22/7.$$

Finally we note that if we find the continued fraction series for $C_{24576} = 3.14159267$, we get the approximants 3/1, 22/7, 333/106, 355/113,... The last number is precisely the approximation discovered by Zu Chongzhi.

## 6. Conclusions

In this article, we have described Archimedes' method of determining the value of $\pi$ by approximating the circumference of a circle of unit diameter by the perimeters of inscribed and circumscribed regular polygons. As we increase n, we get more and more accurate approximations for $\pi$. As we have shown, it turns out that it is easy to find the required perimeters if we restrict ourselves to polygons with number of sides $n=3\times 2^k$, with $k=0,1,2,3...$ that is *n=3, 6, 12, 24, 48...* Archimedes went up to *n=96* and from this he deduced the approximate value *22/7* that is widely used today. Zu Chongzhi went up to *n=24576* and obtained the value *355/113* which is a better approximation of $\pi$. Today we can use a calculator or computer to find the precise numerical value of the perimeter of a polygon, which typically involves finding a lot of square roots. It must have been very difficult in the times of Archimedes and Chongzi to compute these numbers by hand. It would have required a lot of clever thinking to arrive at simple rational approximations such as 22/7 and 355/113 and it is an interesting question as to how exactly they arrived at these results [1].

## 7. Acknowledgements

We thank Ranjini Bandyopadhyay, Deepak Dhar, Rajesh Gopakumar, Anupam Kundu, Joseph Samuel, and Supurna Sinha for their comments and suggestions.

## 8. Suggested Reading